\documentclass[12 pt, a4 paper, twoside]{article}
 \usepackage{amssymb,amsmath,graphicx,enumerate}
 \usepackage{longtable,mathrsfs}
 \usepackage{hyperref}
  \numberwithin{equation}{section}
 \newtheorem{thm}{Theorem}[section]
 \newtheorem{pro}{Propostion}[section]
 \newtheorem{lem}{Lemma}[section]

\newtheorem{re}{Remark}[section]
 \newtheorem{defn}{Definition}[section]

 \allowdisplaybreaks
 \begin{document}
\begin{center}\textbf{LYAPUNOV TYPE INEQUALITY FOR FRACTIONAL DIFFERENTIAL EQUATION WITH k-PRABHAKAR DERIVATIVE}\end{center}
\begin{center}
\textbf{Narayan G. Abuj${^{1}}$ and  Deepak B. Pachpatte${^{2}}$} \\
Department of mathematics,\\ Dr. Babasaheb Ambedkar Marathwada University,\\ Aurangabad-431004 (M.S.) India.\\
{${{^1}}$abujng@gmail.com  and ${{^2}}$pachpatte@gmail.com}
\end{center}
\begin{abstract}
 In this paper, Lyapunov type inequality is establish for fractional boundary value problem involving the k-Prabhakar fractional derivative.
\end{abstract}
${\textbf{\emph{Mathematics Subject Classification}}}: 26D10; 33E12; 34A08$. 
 \section{Introduction}
\paragraph{}Lyapunov type inequality plays an very important role in the study of various properties of solutions of differential and difference equations such as control theory, oscillation theory, disconjugacy, eigenvalue problem etc. {\cite{AML,DCAK,HT,JZ,XY1,XY2,XY3}}. Many authors studied generalization of differential and integral operators such as Prabhakar derivative, Prabhakar integral, k-Riemann-Liouville derivative, k-Riemann-Liouville integral and k-Prabhakar derivative, k-Prabhakar integral. Also obtained Lyapunov type inequalities for some of these differential operators \cite{GD,EA,GCPT,DNA}.
 \paragraph{} In the beginning A. M. Lyapunov considered the boundary value problem (BVP)
  \begin{equation}\label{1}
  \left\{ {\begin{array}{*{20}{l}}{y^{\prime\prime}(t)}+q(t)y(t)= 0,\quad a<t< b,\quad\\
  y(a)=y(b)=0, \end{array}} \right.
  \end{equation}
  has the nontrivial solution then for real and continuous function q(t) he obtained the following well-known result called Lyapunov inequality.  \\
  \begin{equation}\label{2}
    \int_{a}^{b}|q(u)|du > \frac{4}{b-a}.
  \end{equation}
  \paragraph{}Ferreira \cite{RAC1,RAC2} obtained Lypunov inequality for boundary value problem involving Riemann-Liouville and Caputo derivatives.
  Jleli and Samet modified these inequalities of Ferreira \cite{MJ1,MJ2}.
  \paragraph{}Recently S. Eshaghi and A. Ansari \cite{EA} obtained Lyapunov inequality for fractional boundary value problem with Prabhakar derivative, also
  the authors D. B Pachpatte and et.al \cite{DNA} developed Lyapunov type inequality for hybrid fractional boundary value problem involving Prabhakar derivative.
 \paragraph{} In this paper, we consider the following fractional boundary value problem involving k-Prabhakar derivative
 \begin{equation}\label{3}
 \left\{ {\begin{array}{*{20}{l}}\textbf(_{k}{D}_{\rho, \beta, \omega, a+}^{\gamma}y)(t)+q(t)y(t)= 0,\quad a<t<b,  \\
 y(a)=y(b)=0, \end{array}}\right.
 \end{equation}
 Where $y\in C[a,b]$, $_{k}{D}_{\rho, \beta, \omega, a+}^{\gamma}$ is k-Prabhakar differntial operator of order $\beta\in(1,2]$, $k\in\mathbb{R}^{+} and \quad \rho, \gamma, \omega \in \mathbb{C}$.
 We obtained Greens function of the fractional boundary value problem \eqref{3} in terms of k-Mittag-Leffler function. Also we state and prove some properties of Green function and establish the lyapunov inequality for the fractional boundary value problem \eqref{3}.
\section{Preliminaries}
 In this section, we collect some basic definitions and lemmas that will be important to us in the sequel.
 \begin{defn}\cite{DC}
  The k-Mittag-Leffler function is denoted by $E_{k,\alpha,\beta}^{\gamma}(z)$ and is defined as
 \begin{equation}\label{4}
 E_{k,\alpha,\beta}^{\gamma}(z)= \sum_{n=0}^\infty\frac{(\gamma)_{n,k}\,z^n}{\Gamma_{k}(\alpha n +\beta) n!}
 \end{equation}
 where $ k\in \mathbb{R}^{+}, \alpha, \beta, \gamma\in \mathbb{C}$, $\Re(\alpha)>0$, $\Re(\beta)>0$; $\Gamma_{k}(x)$ is the k-Gamma function and $(\gamma)_{n,k}=\frac{\Gamma_{k}(\gamma+nk)}{\Gamma_{k}(\gamma)}$ is the pochhammer k-symbol. \\
 \end{defn}
\begin{defn}\cite{GD}
 Let $\alpha, \beta ,\omega,\gamma\in\mathbb{C}$, $k\in \mathbb{R^{+}}$, $\Re(\alpha)>0$, $\Re(\beta)>0$ and $\phi\in L^{1}([0,b])$,$(0<x<b\leq\infty)$. The k-Prabhakar integral operator involving k-Mittag-Leffler function is defined as
\begin{align}\label{5}
    (_{k}\textbf{P}_{\alpha,\beta,\omega}\phi)(x)&=\int_{0}^{x}\frac{(x-t)^{\frac{\beta}{k}-1}}{k}E_{k,\alpha,\beta}^{\gamma}[\omega(x-t)^{\frac{\alpha}{k}}]\phi(t)dt,\quad(x>0)\\
   & = (_{k}\mathcal{E}_{\alpha,\beta,\omega}^{\gamma}\ast f)(x)\\
   where \nonumber
  \end{align}
  \begin{equation}\label{6}
   \,\,\,\,\,\,{_{k}\mathcal{E}_{\alpha,\beta,\omega}^{\gamma}}(t)= \left\{ {\begin{array}{*{20}{l}}{\frac{t^{\frac{\beta}{k}-1}}{k}}E_{k,\alpha,\beta}^{\gamma}(\omega t^{\frac{\alpha}{k}}),\quad\quad t>0; \\
  0, \qquad \qquad \qquad\qquad t\leq 0. \end{array}} \right.\\
  \end{equation}
 \end{defn}
 \begin{defn}\cite{GD}
   Let $k\in \mathbb{R^{+}}$, $\rho,\beta,\gamma,\omega\in\mathbb{C}$,\,$\Re(\alpha)>0$,\,$\Re(\beta)>0$,\\
   $m=[\frac{\beta}{k}]+1$, $f\in L^{1}([0,b])$. The k-Prabhakar derivative is defined as
 \begin{equation}\label{a1}
 _{k}\textbf{D}_{\rho,\beta,\omega,}^{\gamma}f(x)=\bigg(\frac{d}{dx}\bigg)^{m}k^{m}\,_{k}\textbf{P}_{\rho, mk-\beta,\omega}^{-\gamma}f(x).
 \end{equation}
 \end{defn}
  \begin{lem}\label{L1}\cite{GD}  Let $\alpha,\beta,\omega,\gamma\in\mathbb{C}$,\,$ k\in \mathbb{R^{+}}$,\,$\Re(\alpha)>0; \Re(\beta)>0$, $\phi\in L^{1}(\mathbb{R}_{0}^{+})$ and ${|\omega k(ks)^{\frac{-\rho}{k}}|<1}$ then
 \begin{align}\label{c1}
 \mathscr{L}\{(_{k}\textbf{P}_{\rho,\beta,\omega,}^{\gamma}\phi)(x)\}(s)&= \mathscr{L}\{_{k}\varepsilon_{\rho,\beta,\omega}^{\gamma}(t)\}(s)\mathscr{L}\{\phi\}(s)\nonumber\\
 &=(ks)^{\frac{-\beta}{k}}(1-\omega k(ks)^{\frac{-\rho}{k}})^{\frac{-\gamma}{k}} \mathscr{L}\{\phi\}(s).
 \end{align}
 \end{lem}
 \begin{lem}\label{L2}\cite{PKD}
 The Laplace transform of k-Prabhakar derivative \eqref{a1} is
 \begin{equation}\label{10}\begin{split}
 \mathscr{L}\{_{k}\textbf{D}_{\rho,\beta,\omega,}^{\gamma}f(x)\}=&(ks)^{\frac{\beta}{k}}(1-\omega k(ks)^{\frac{-\rho}{k}})^{\frac{\gamma}{k}}F(s)\\
 &-\sum_{n=0}^{m-1}k^{n+1}s^{n}\Bigg({_{k}\textbf{D}_{\rho,\beta-(n+1)k,\omega}^{\gamma}f(0^{+})}\Bigg).
 \end{split}\end{equation}
 For the case $[\frac{\beta}{k}]+1= m=1$,
 \begin{equation}\label{9}
 \mathscr{L}\{_{k}\textbf{D}_{\rho,\beta,\omega,}^{\gamma}y(x)\}=(ks)^{\frac{-\beta}{k}}(1-\omega k(ks)^{\frac{-\rho}{k}})^{\frac{\gamma}{k}} \mathscr{L}\{y(x)\}(s)-k(_{k}\textbf{P}_{\rho,k-\beta,\omega,}^{-\gamma}y)(0)\nonumber
 \end{equation}
 with ${|\omega k(ks)^{\frac{-\rho}{k}}|<1}$.
 \end{lem}
 \section{Main Results}
 In this section, we shall establish our Lyapunov type inequality with the help of following propositions.
 \begin{pro}\label{P1}
 If $f(x)\in C(a,b)\cap L(a,b)$; then ${_{k}\textbf{D}_{\rho,\beta,\omega,a+}^{\gamma}}{_{k}\textbf{P}_{\rho,\beta,\omega,a+}^{\gamma}}f(x)=f(x)$
 and if $f(x),{_{k}\textbf{D}_{\rho,\beta,\omega,a+}^{\gamma}}f(x) \in C(a,b)\cap L(a,b)$, then for $c_{j}\in{\mathbb{R}}$ and $m-1<\beta\leq m$, we have
 \begin{align}\label{11}
 \hspace{-.5cm}{_{k}\textbf{P}_{\rho,\beta,\omega,a+}^{\gamma}}{_{k}\textbf{D}_{\rho,\beta,\omega,a+}^{\gamma}}f(x) = & f(x)+ c_{0}(x-a)^{\frac{\beta}{k}-1}{E_{k,\rho,\beta}^{\gamma}(\omega(x-a)^{\frac{\rho}{k}})} \nonumber\\
   &+  c_{1}(x-a)^{\frac{\beta}{k}-2}{E_{k,\rho,\beta-k}^{\gamma}(\omega(x-a)^{\frac{\rho}{k}})} \nonumber \\
   &+c_{2}(x-a)^{\frac{\beta}{k}-3}{E_{k,\rho,\beta-2k}^{\gamma}(\omega(x-a)^{\frac{\rho}{k}})}+ ... \nonumber\\
   &+c_{m-1}(x-a)^{\frac{\beta}{k}-m}{E_{k,\rho,\beta-(m-1)k}^{\gamma}(\omega(x-a)^{\frac{\rho}{k}})}
\end{align}
\end{pro}
\textbf{Proof}. With the advantage of relation \eqref{c1} and \eqref{10} Laplace transform of L.H.S of \eqref{11} for $a=0$ is
\begin{align*}\label{12}
\mathscr{L}\Bigg\{{_{k}\textbf{P}_{\rho,\beta,\omega,0+}^{\gamma}}{_{k}\textbf{D}_{\rho,\beta,\omega,0+}^{\gamma}}f(x); s\Bigg\}= F(s)&- k\sum_{n=0}^{m-1}(ks)^{\frac{-\beta+nk}{k}}(1-\omega k(ks)^{\frac{-\rho}{k}})^{\frac{-\gamma}{k}}\\
&\times({_{k}\textbf{D}_{\rho,\beta-(n+1)k,\omega,}^{\gamma}}f)(0^{+})\\
\end{align*} Inverting Laplace of above equation for $a\neq0$ gives the desired proof.\\
\begin{re}
 Note that, for $\beta=\mu,k=1,$ equation \eqref{11} coincides with (\cite{EA} equation (27)).
\end{re}
\begin{pro}\label{P2}
Let $ k\in \mathbb{R^{+}}$,$\rho,\beta,\gamma,\omega \in\mathbb{C}$, $\Re(\alpha)>0; \,\Re(\beta)>0$ then for any $j\in\mathbb{N}$ we have
\begin{equation}\label{13}
  \frac{d^{j}}{dx^{j}}[x^{\frac{\beta}{k}-1}{E_{k,\rho,\beta}^{\gamma}(\omega x^{\frac{\rho}{k}})}]=\frac{x^{\frac{\beta}{k}-(j+1)}}{k^j}E_{k,\rho,\beta-jk}^{\gamma}(\omega x^{\frac{\rho}{k}})
\end{equation}
\end{pro}
\textbf{Proof}. We prove this result by using property of k-gamma function \cite{CGK} as,
\begin{align*}\label{14}
 {\frac{d}{dx}}[x^{\frac{\beta}{k}-1}{E_{k,\rho,\beta}^{\gamma}(\omega x^{\frac{\rho}{k}})}]&= {\frac{d}{dx}}\Bigg(x^{\frac{\beta}{k}-1}\sum_{n=0}^\infty\frac{(\gamma)_{n,k}\, \omega^n x^{\frac{\rho n}{k}}}{\Gamma_{k}(\rho n +\beta) n!}\Bigg)\nonumber\\
 &=\sum_{n=0}^\infty\frac{(\gamma)_{n,k}\, \omega^n }{\Gamma_{k}(\rho n +\beta) n!}{\frac{d}{dx}x^{\frac{\rho n+\beta-k}{k}}}\nonumber\\
 &= \sum_{n=0}^\infty\frac{(\gamma)_{n,k}\, \omega^n }{k\Gamma_{k}(\rho n +\beta) n!}  (\rho n+\beta-k)x^{\frac{\rho n+\beta-k}{k}-1}\nonumber\\
 &= \sum_{n=0}^\infty\frac{(\gamma)_{n,k}\, \omega^n }{k\Gamma_{k}(\rho n +\beta-k+k) n!}  (\rho n+\beta-k)x^{\frac{\rho n+\beta-2k}{k}}\nonumber\\
  &= x^{\frac{\beta-2k}{k}}\sum_{n=0}^\infty\frac{(\gamma)_{n,k} (\omega^n x^{\frac{\rho n}{k}})}{k\Gamma_{k}(\rho n +\beta-k) n!}\nonumber \\
  &= \frac{x^{\frac{\beta}{k}-2}}{k}{E_{k,\rho,\beta-k}^{\gamma}(\omega x^{\frac{\rho}{k}})}\\
\end{align*}
Similarly, we find the second derivative  \\
\begin{align*}
\frac{d^{2}}{dx^{2}}[x^{\frac{\beta}{k}-1}{E_{k,\rho,\beta}^{\gamma}(\omega x^{\frac{\rho}{k}})}]=
\frac{x^{\frac{\beta}{k}-3}}{k^2}{E_{k,\rho,\beta-2k}^{\gamma}(\omega x^{\frac{\rho}{k}})}
\end{align*}
continuing this process j-times we have the desired result.
\begin{re}
Note that, for $j=n,\beta=\mu$ and $k=1,$ equation \eqref{13} coincides with (\cite{EA}, equation (14)).
\end{re}
\begin{thm}\label{thm1}
Let $1<\beta\leq2$,\,$\gamma,\rho,\omega\in\mathbb{R^{+}}$,\,$y\in\mathbf{C}[a,b]\cap L[a,b]$,\,then the fractional boundary value problem
\end{thm}
\begin{equation}\label{15}
 \left\{ {\begin{array}{*{20}{l}}\textbf(_{k}{D}_{\rho, \beta, \omega, a+}^{\gamma}y)(t)+q(t)y(t)= 0,\quad a<t<b,  \\
 y(a)=y(b)=0. \end{array}}\right.
 \end{equation}
 is equivalent to the following integral equation
\begin{equation}\label{16}
y(t)=\int_{a}^{b}G(t,u)q(u)y(u)du,
\end{equation}
where the G(t,u) is the Greens function and is given by
\begin{equation}\label{17}
G(t,u)= \left\{
{\begin{array}{*{20}{l}}\frac{(t-a)^{{\frac{\beta}{k}}-1}E_{k,\rho,\beta}^{\gamma}(\omega(t-a)^{\frac{\rho}{k}})}{(b-a)^{{\frac{\beta}{k}}-1}E_{k,\rho,\beta}^{\gamma}(\omega(b-a)^{\frac{\rho}{k}})}\frac{(b-u)^{{\frac{\beta}{k}}-1}}{k}{E_{k,\rho,\beta}^{\gamma}(\omega(b-u)^{\frac{\rho}{k}})}\\
 -\frac{(t-u)^{{\frac{\beta}{k}}-1}}{k}E_{k,\rho,\beta}^{\gamma}(\omega(t-u)^{\frac{\rho}{k}}),\quad\quad\quad\quad\quad\quad\quad\quad \quad a\leq u\leq t\leq b, \\
\frac{(t-a)^{{\frac{\beta}{k}}-1}E_{k,\rho,\beta}^{\gamma}(\omega(t-a)^{\frac{\rho}{k}})}{(b-a)^{{\frac{\beta}{k}}-1}E_{k,\rho,\beta}^{\gamma}(\omega(b-a)^{\frac{\rho}{k}})}\frac{(b-u)^{{\frac{\beta}{k}}-1}}{k}{E_{k,\rho,\beta}^{\gamma}(\omega(b-u)^{\frac{\rho}{k}})},\,\ \,a\leq t\leq u\leq b. \end{array}} \right.
\end{equation}
\textbf{Proof}. Applying the k-Prabhakar integral operator  $_{k}\textbf{P}_{\rho,\beta,\omega,a+}^{\gamma}$ on the fractional differential equation \eqref{15} and using proposition \ref{P1} we find real constant as follows
\begin{align}\label{18}
\hspace{-0.5cm} y(t)&= -\int_{a}^{t}\frac{(t-u)^{\frac{\beta}{k}-1}}{k}E_{k,\rho,\beta}^{\gamma}[\omega(t-u)^{\frac{\rho}{k}}]q(u)y(u)du\nonumber\\ &+c_{0}(t-a)^{\frac{\beta}{k}-1}E_{k,\rho,\beta}^{\gamma}[\omega(t-a)^{\frac{\rho}{k}}]
  +c_{1}(t-a)^{\frac{\beta}{k}-2}E_{k,\rho,\beta-k}^{\gamma}[\omega(t-a)^{\frac{\rho}{k}}].
\end{align}
By employing the boundary conditions yields $c_{1}=0$ and $$c_{0}=\frac{1}{{(b-a)^{\frac{\beta}{k}-1}}E_{k,\rho,\beta}^{\gamma}(\omega(b-a)^{\frac{\rho}{k}})}
\int_{a}^{b}\frac{(b-u)^{\frac{\beta}{k}-1}}{k}E_{k,\rho,\beta}^{\gamma}[\omega(b-u)^{\frac{\rho}{k}}]q(u)y(u)du.$$
Substituting these values of real constants in equation \eqref{18} we have the unique solution of \eqref{15} as follows
\newpage
\begin{align*}
  y(t)&=\int_{a}^{t}\Bigg[\frac{(t-a)^{{\frac{\beta}{k}}-1}E_{k,\rho,\beta}^{\gamma}(\omega(t-a)^{\frac{\rho}{k}})}{(b-a)^{{\frac{\beta}{k}}-1}E_{k,\rho,\beta}^{\gamma}(\omega(b-a)^{\frac{\rho}{k}})}\frac{(b-u)^{{\frac{\beta}{k}}-1}}{k}{E_{k,\rho,\beta}^{\gamma}(\omega(b-u)^{\frac{\rho}{k}})}\\
 &-\frac{(t-u)^{{\frac{\beta}{k}}-1}}{k}E_{k,\rho,\beta}^{\gamma}(\omega(t-u)^{\frac{\rho}{k}})\Bigg]q(u)y(u)du\\
 &+\int_{t}^{b}\Bigg[\frac{(t-a)^{{\frac{\beta}{k}}-1}E_{k,\rho,\beta}^{\gamma}(\omega(t-a)^{\frac{\rho}{k}})(b-u)^{{\frac{\beta}{k}}-1}}{k(b-a)^{{\frac{\beta}{k}}-1}E_{k,\rho,\beta}^{\gamma}(\omega(b-a)^{\frac{\rho}{k}})}{E_{k,\rho,\beta}^{\gamma}(\omega(b-u)^{\frac{\rho}{k}})}\Bigg]q(u)y(u)du\\
  y(t)&=\int_{a}^{b}G(t,u)q(u)y(u)du.
\end{align*}
Where G(t,u) is Greens function given by \eqref{17}.
\begin{thm}\label{thm2}
The Greens function defined by \eqref{17} holds the following properties:\\
\textbf{(a)}\quad $G(t,u)\geq0$,\quad$\forall$\,$a\leq t,u \leq b$.\\
\textbf{(b)}\quad $\mathop {\max }\limits_{t\in[a,b]}G(t,u)= G(u,u)$, for $ u\in[a,b]$.\\
\textbf{(c)}\quad The maximum of the function $G(u,u)$ occurs at point $u=\frac{a+b}{2}$ and has the maximum value is \\
\begin{equation}\label{19}
\mathop{\max}\limits_{u\in[a,b]}G(u,u)= \bigg{(}\frac{b-a}{4}\bigg{)}^{{\frac{\beta}{k}}-1}\frac{E_{k,\rho,\beta}^{\gamma}(\omega(\frac{b-a}{2})^{\frac{\rho}{k}})
E_{k,\rho,\beta}^{\gamma}(\omega(\frac{b-a}{2})^{\frac{\rho}{k}})}{k.E_{k,\rho,\beta}^{\gamma}(\omega(b-a)^{\frac{\rho}{k}})}.
\end{equation}
\end{thm}
\textbf{Proof}. We Prove this theorem by setting two function as follows\\
\begin{align*}
 g _{1}(t,u)&= \frac{(t-a)^{{\frac{\beta}{k}}-1}E_{k,\rho,\beta}^{\gamma}(\omega(t-a)^{\frac{\rho}{k}})}{(b-a)^{{\frac{\beta}{k}}-1}E_{k,\rho,\beta}^{\gamma}
 (\omega(b-a)^{\frac{\rho}{k}})}\frac{(b-u)^{{\frac{\beta}{k}}-1}}{k}{E_{k,\rho,\beta}^{\gamma}(\omega(b-u)^{\frac{\rho}{k}})}\\ &-\frac{(t-u)^{{\frac{\beta}{k}}-1}}{k}E_{k,\rho,\beta}^{\gamma}(\omega(t-u)^{\frac{\rho}{k}}),\quad a\leq u\leq t\leq b,\\
 and \\
 g _{2}(t,u)&= \frac{(t-a)^{{\frac{\beta}{k}}-1}E_{k,\rho,\beta}^{\gamma}(\omega(t-a)^{\frac{\rho}{k}})}{(b-a)^{{\frac{\beta}{k}}-1}E_{k,\rho,\beta}^{\gamma}
(\omega(b-a)^{\frac{\rho}{k}})}\frac{(b-u)^{{\frac{\beta}{k}}-1}}{k}{E_{k,\rho,\beta}^{\gamma}(\omega(b-u)^{\frac{\rho}{k}})}, \,\,a\leq t\leq u\leq b\\
\end{align*}
It is obvious that $ g_{2}(t,u)\geq 0$. So to prove (a), we have to show that $g_{1}(t,u)\geq 0$ or it is equivalent to prove that\\
$\frac{(t-a)^{{\frac{\beta}{k}}-1}E_{k,\rho,\beta}^{\gamma}(\omega(t-a)^{\frac{\rho}{k}})}{(b-a)^{{\frac{\beta}{k}}-1}E_{k,\rho,\beta}^{\gamma}
 (\omega(b-a)^{\frac{\rho}{k}})}\frac{(b-u)^{{\frac{\beta}{k}}-1}}{k}{E_{k,\rho,\beta}^{\gamma}(\omega(b-u)^{\frac{\rho}{k}})} \geq\frac{(t-u)^{{\frac{\beta}{k}}-1}}{k}E_{k,\rho,\beta}^{\gamma}(\omega(t-u)^{\frac{\rho}{k}})$.\\
 \newpage
Therefore it is sufficient to prove that \\
$(i)$\qquad$\frac{(t-a)^{\frac{\beta}{k}-1}(b-u)^{\frac{\beta}{k}-1}}{k(b-a)^{\frac{\beta}{k}-1}}\geq \frac{(t-u)^{\frac{\beta}{k}-1}}{k}$ or $\frac{(t-a)^{\frac{\beta}{k}-1}(b-u)^{\frac{\beta}{k}-1}}{(b-a)^{\frac{\beta}{k}-1}}\geq {(t-u)^{\frac{\beta}{k}-1}}$,\\
$(ii)$\qquad$\frac{E_{k,\rho,\beta}^{\gamma}(\omega(t-a)^{\frac{\rho}{k}})E_{k,\rho,\beta}^{\gamma}(\omega(b-u)^{\frac{\rho}{k}})}{E_{k,\rho,\beta}^{\gamma}(\omega(b-a)^{\frac{\rho}{k}})}\geq E_{k,\rho,\beta}^{\gamma}(\omega(t-u)^{\frac{\rho}{k}})$.\\
For proof (i) we proceed as \\
$\frac{(t-a)^{\frac{\beta}{k}-1}(b-u)^{\frac{\beta}{k}-1}}{(b-a)^{\frac{\beta}{k}-1}}\geq (t-u)^{\frac{\beta}{k}-1}$\\
$\Leftrightarrow$\,$\frac{(t-a)^{\frac{\beta}{k}-1}}{(b-a)^{\frac{\beta}{k}-1}}(b-u)^{\frac{\beta}{k}-1}\geq \frac{(t-a)^{\frac{\beta}{k}-1}}{(b-a)^{\frac{\beta}{k}-1}}\Bigg[b-\Bigg(a+\frac{(u-a)(b-a)}{t-a}\Bigg)\Bigg]^{\frac{\beta}{k}-1}$\\
$\Leftrightarrow$\,$a+\frac{(u-a)(b-a)}{t-a}\geq u$\\
$\Leftrightarrow$\,$\frac{at-ab+bu-au}{t-a}\geq u$\\
$\Leftrightarrow$\,$a(t-b)+u(b-t)\geq 0$\\
$\Leftrightarrow$\,$u(b-t)\geq a(b-t)$\\
$\Leftrightarrow$\,$u\geq a$.\\
According to inequality $(t-a)(b-u)\geq (b-a)(t-u)$ and Taylor series expansion of the generalized k-Mittag Leffler function $E_{k,\rho,\beta}^{\gamma}(z)$, for $1<\beta\leq2,$ $\gamma,\rho,\omega,z\in \mathbb{R^{+}}$, Hence the proof of (ii) is complete.\\
 \textbf{Proof of (b)}: We prove this as follows \\
 Differentiate $g_{1}(t,u)$ with respect to t keeping u fixed and apply the proposition \ref{P2} for $j=1$ we have
 \begin{align*}
 g^{\prime}_{1}(t,u)&=\frac{(b-u)^{\frac{\beta}{k}-1}E_{k,\rho,\beta}^{\gamma}(\omega(b-u)^{\frac{\rho}{k}})}
 {k(b-a)^{\frac{\beta}{k}-1}E_{k,\rho,\beta}^{\gamma}(\omega(b-a)^{\frac{\rho}{k}})}\frac{d}{dt}(t-a)^{\frac{\beta}{k}-1} E_{k,\rho,\beta}^{\gamma}(\omega(t-a)^{\frac{\rho}{k}})\\
 &-\frac{1}{k}\frac{d}{dt}(t-u)^{\frac{\beta}{k}-1}E_{k,\rho,\beta}^{\gamma}(\omega(t-u)^{\frac{\rho}{k}})\\ &=\frac{(b-u)^{\frac{\beta}{k}-1}E_{k,\rho,\beta}^{\gamma}(\omega(b-u)^{\frac{\rho}{k}})}
 {k^{2}(b-a)^{\frac{\beta}{k}-1}E_{k,\rho,\beta}^{\gamma}(\omega(b-a)^{\frac{\rho}{k}})}(t-a)^{\frac{\beta}{k}-2}E_{k,\rho,\beta-k}^{\gamma}(\omega(t-a)^{\frac{\rho}{k}})\\
 &-\frac{1}{k^{2}}(t-u)^{\frac{\beta}{k}-2}E_{k,\rho,\beta-k}^{\gamma}(\omega(t-u)^{\frac{\rho}{k}})\\
 u\geq a\Rightarrow\\
 g^{\prime}_{1}(t,u)&\leq\frac{(b-a)^{\frac{\beta}{k}-1}E_{k,\rho,\beta}^{\gamma}(\omega(b-a)^{\frac{\rho}{k}})}
 {k^{2}(b-a)^{\frac{\beta}{k}-1}E_{k,\rho,\beta}^{\gamma}(\omega(b-a)^{\frac{\rho}{k}})}(t-a)^{\frac{\beta}{k}-2}E_{k,\rho,\beta-k}^{\gamma}(\omega(t-a)^{\frac{\rho}{k}})\\
 &-\frac{1}{k^{2}}(t-a)^{\frac{\beta}{k}-2}E_{k,\rho,\beta-k}^{\gamma}(\omega(t-a)^{\frac{\rho}{k}})\\
 &\leq 0\\
 \end{align*}
 that yields $g_{1}(t,u)$ is a decreasing function of t. Similarly by differentiating $g_{2}(t,u)$ with respect to t for every fixed u. From this we conclude that $g_{2}(t,u)$ is increasing function.
 Therefore, the maximum  of the function $G(t,u)$ with respect to t is the value G(u,u).
 Finally, we set the function $h(u)$ for $u\in[a,b]$ as follows
 \begin{align}\label{20}
   h(u) = G(u,u)= \frac{(u-a)^{{\frac{\beta}{k}}-1}E_{k,\rho,\beta}^{\gamma}(\omega(u-a)^{\frac{\rho}{k}})}{(b-a)^{{\frac{\beta}{k}}-1}E_{k,\rho,\beta}^{\gamma}
(\omega(b-a)^{\frac{\rho}{k}})}\frac{(b-u)^{{\frac{\beta}{k}}-1}}{k}{E_{k,\rho,\beta}^{\gamma}(\omega(b-u)^{\frac{\rho}{k}})} \nonumber\\
\end{align}
By using the equation \eqref{4} and \eqref{13} we have
\begin{align*}
 h^{\prime}(u)&=\frac{1}{k(b-a)^{{\frac{\beta}{k}}-1E_{k,\rho,\beta}^{\gamma}(\omega(b-a)^{\frac{\rho}{k}})}}\frac{d}{du}
[(u-a)^{\frac{\beta}{k}-1}E_{k,\rho,\beta}^{\gamma}(\omega(u-a)^{\frac{\rho}{k}})\\
&\hspace{2cm}\times(b-u)^{\frac{\beta}{k}-1}E_{k,\rho,\beta}^{\gamma}(\omega(b-u)^{\frac{\rho}{k}})]\\
&=\frac{((u-a)(b-u))^{\frac{\beta}{k}-2}}{k^{2}(b-a)^{{\frac{\beta}{k}}-1E_{k,\rho,\beta}^{\gamma}(\omega(b-a)^{\frac{\rho}{k}})}}
\Bigg[(b-u)E_{k,\rho,\beta-k}^{\gamma}(\omega(u-a)^{\frac{\rho}{k}})\\
&\times E_{k,\rho,\beta}^{\gamma}(\omega(b-u)^{\frac{\rho}{k}})-(u-a)E_{k,\rho,\beta-k}^{\gamma}(\omega(b-u)^{\frac{\rho}{k}})
E_{k,\rho,\beta}^{\gamma}(\omega(u-a)^{\frac{\rho}{k}})\Bigg]\\
&=\frac{((u-a)(b-u))^{\frac{\beta}{k}-2}}{k^{2}(b-a)^{{\frac{\beta}{k}}-1E_{k,\rho,\beta}^{\gamma}(\omega(b-a)^{\frac{\rho}{k}})}}
\Bigg[(b-u)\sum_{n=0}^\infty\frac{(\gamma)_{n,k} (\omega(u-a))^{\frac{\rho n}{k}})}{\Gamma_{k}(\rho n +\beta-k) n!}\\
&\hspace{1cm}\times\sum_{n=0}^\infty\frac{(\gamma)_{n,k} (\omega(b-u))^{\frac{\rho n}{k}})}{\Gamma_{k}(\rho n +\beta) n!}
-(u-a)\sum_{n=0}^\infty\frac{(\gamma)_{n,k} (\omega(b-u))^{\frac{\rho n}{k}})}{\Gamma_{k}(\rho n +\beta-k) n!}\\
&\hspace{1cm}\times\sum_{n=0}^\infty\frac{(\gamma)_{n,k} (\omega(u-a))^{\frac{\rho n}{k}})}{\Gamma_{k}(\rho n +\beta) n!}\Bigg]
\end{align*}
By solving $h^{\prime}(u)=0\Rightarrow u=\frac{a+b}{2}$, also we observe that $h^{\prime}(u)>0$ on $(a,\frac{a+b}{2})$ and $h^{\prime}(u)<0$ on $(\frac{a+b}{2},b)$. Hence h(u) has maximum at point $u=\frac{a+b}{2}$.\\
\textbf{Proof of (c)}: By substituting $u=\frac{a+b}{2}$ in \eqref{20} it gives the maximum of G(u,u) as follows
\begin{equation}\label{21}
\mathop{\max}\limits_{u\in[a,b]}G(u,u)=
\Bigg(\frac{b-a}{4}\Bigg)^{\frac{\beta}{k}-1}\frac{E_{k,\rho,\beta}^{\gamma}(\omega{(\frac{b-a}{2})^{\frac{\rho}{k}}})
E_{k,\rho,\beta}^{\gamma}(\omega{(\frac{b-a}{2})^{\frac{\rho}{k}}})}{k\,E_{k,\rho,\beta}^{\gamma}(\omega(b-a)^{\frac{\rho}{k}})}.
\end{equation}
\begin{thm}\label{thm3} Let $\mathscr{B}=C[a,b]$ be the Banach space equipped with norm $ \|y\|= \mathop{\sup}\limits_{t\in[a,b]}|y(t)|$ and a nontrivial continuous solution of the fractional boundary value problem \eqref{3} exists then
\begin{equation}\label{22}
\int_{a}^{b}|q(u)|du\geq\Bigg(\frac{4}{b-a}\Bigg)^{\frac{\beta}{k}-1}\frac{k\,E_{k,\rho,\beta}^{\gamma}(\omega(b-a)^{\frac{\rho}{k}})}
{E_{k,\rho,\beta}^{\gamma}(\omega{(\frac{b-a}{2})^{\frac{\rho}{k}}})E_{k,\rho,\beta}^{\gamma}(\omega{(\frac{b-a}{2})^{\frac{\rho}{k}}})}\\
\end{equation}
where q(t) is real and continuous function.
\end{thm}
\textbf{Proof}. According to theorem \ref{thm1} a solution of the above fractional boundary value problem \eqref{3} satisfies the integral equation\\
\begin{align*}
y(t)=\int_{a}^{b}G(t,u)q(u)y(u)du\\
\end{align*}
which by applying the indicated norm on both side of it and using the second and third properties of theorem \eqref{thm2} we get  the desired inequality as follows
\begin{align*}
\|y(t)\|&\leq\mathop{\max}\limits_{t\in[a,b]}\int_{a}^{b}|G(t,u)q(u)|du\|y(t)\|\\
1&\leq\mathop{\max}\limits_{t\in[a,b]}\int_{a}^{b}|G(t,u)q(u)|du\\
1&\leq\Bigg(\frac{b-a}{4}\Bigg)^{\frac{\beta}{k}-1}\frac{E_{k,\rho,\beta}^{\gamma}(\omega{(\frac{b-a}{2})^{\frac{\rho}{k}}})
E_{k,\rho,\beta}^{\gamma}(\omega{(\frac{b-a}{2})^{\frac{\rho}{k}}})}{k\,E_{k,\rho,\beta}^{\gamma}(\omega(b-a)^{\frac{\rho}{k}})}\int_{a}^{b}|q(u)|du\\
\int_{a}^{b}|q(u)|du&\geq\Bigg(\frac{4}{b-a}\Bigg)^{\frac{\beta}{k}-1}\frac{k\,E_{k,\rho,\beta}^{\gamma}(\omega(b-a)^{\frac{\rho}{k}})}
{E_{k,\rho,\beta}^{\gamma}(\omega{(\frac{b-a}{2})^{\frac{\rho}{k}}})E_{k,\rho,\beta}^{\gamma}(\omega{(\frac{b-a}{2})^{\frac{\rho}{k}}})}.\\
\end{align*}
\section{Conclusion}
In this paper, we obtained more general results than in \cite{EA}. The results in \cite{EA} can be obtained for particular values of k and $\beta$ as $k=1$ and $\beta=\mu$ in Greens function in theorem \ref{thm1}. and Lyapunov inequality in theorem \ref{thm3}.

\end{document}